\documentclass[11pt]{amsart}
\usepackage{latexsym,amssymb,amsmath}
\def\hd{,...,}
\def\trank{\text{rank}}
\def\tdim{\text{dim}}
\def\fa{{\mathfrak a}}

\def\CC{\mathbb C}
\def\RR{\mathbb R}
\def\HH{\mathbb H}
\def\AA{{\mathbb A}}
\def\BB{{\mathbb B}}
\def\OO{\mathbb O}

\def\ZZ{\mathbb Z}

\def\11{\mathbf 1}
\def\PP{\mathbb P}

\def\fh{{\mathfrak h}}

\def\fd{{\mathfrak d}}
\def\fsl{{\mathfrak {sl}}}
\def\fsp{{\mathfrak {sp}}}

\def\fso{{\mathfrak {so}}}
\def\fe{{\mathfrak e}}
\def\ff{{\mathfrak f}}

\def\fg{{\mathfrak g}}

\def\ft{{\mathfrak t}}

\def\l{\lambda}
\def\a{\alpha}

\def\o{\omega}

\def\b{\beta}
\def\g{\gamma}
\def\s{\sigma}

\def\d{\delta}
\def\th{\theta}
\def\m{\mu}
\def\up#1{{}^{({#1})}}
\def\e{\varepsilon}
\def\ot{{\mathord{\,\otimes }\,}}
\def\op{{\mathord{\,\oplus }\,}}

\def\ra{{\mathord{\;\rightarrow\;}}}

\def\dim{{\rm dim}\;}
\def\La{\Lambda}

\newcommand\exam{{\medskip\noindent {\em Example}.}\hspace{2mm}}

\newtheorem{theo}{Theorem}[section]
\newtheorem{coro}[theo]{Corollary}
\newtheorem{lemm}[theo]{Lemma}
\newtheorem{prop}[theo]{Proposition}

\begin{document}

\title{Triality, exceptional Lie algebras \linebreak
and Deligne dimension formulas}
\author{J.M. Landsberg and L. Manivel}
\date{June 2001}

\begin{abstract} We give a computer free proof of the Deligne, Cohen
and deMan formulas for the dimensions of the irreducible $\fg$-modules
appearing in $\fg^{\ot k}$, $k\leq 4$, where $\fg$ ranges
over  the exceptional complex
simple Lie algebras. We give additional dimension formulas for the
exceptional series, as well as uniform dimension formulas for other
representations distinguished by Freudenthal along the rows of 
his magic chart. Our proofs use the triality model of the magic 
square which we review and present a simplified proof of its
validity. We conclude
with some general remarks about obtaining \lq\lq series\rq\rq\ of
Lie algebras in the spirit of Deligne and Vogel.\end{abstract} 

\maketitle

\section{Introduction}

The goal of this paper is to give a partial  explanation to some 
astonishing observations made by Deligne about the exceptional
complex simple Lie algebras \cite{del}. 
Deligne, following a remark of Vogel, noticed that the tensor powers
$\fg^{\ot k}$ for $\fg$ an  exceptional complex simple  Lie algebra, 
decomposed uniformly into irreducible $\fg$-modules when $k\leq 4$. 
Parametrizing the exceptional series $\fa_1,\fa_2,\fg_2,\fd_4,
\ff_4,\fe_6,\fe_7,\fe_8$ by the inverse Coxeter number $\l$, 
he, together with Cohen and de Man, gave the dimensions of the corresponding
irreducible modules in terms of rational functions of $\l$. 
These rational functions, computed by LiE \cite{LiE}, had 
the ``miraculous\rq\rq\
property that both the numerators and denominators were products of 
{\it linear} functions of $\l$.

Inspired by work of Freudenthal and Tits, we thought it might be interesting
to parametrize the exceptional series by $a=\dim_{\CC}\AA$,
where $\AA$ is respectively the complexification of $0,\RR,\CC,\HH,\OO$
for the last five algebras in the exceptional series (so $a=0,1,2,4,8$).
A first indication that this might be fruitful was the simple 
relation  $\l=-\frac 2{a+2}$.
 The parameter $a$ simplified the Deligne
dimension formula because every time a power of $\l$ appears in the
denominator (which is always), its contribution to the degree of the
denominator is erased upon the change of variable, so that
using $a$,  the denominators
have lower degree and the numerators the same degree. 

The presence of only linear 
forms in the Deligne dimension formulas also suggests one should attempt
to apply the  Weyl dimension 
formula in a uniform way, which is what we have done.

To do this, we found a suitable variant of the  Vinberg 
construction of the exceptional Lie algebras in terms of normed 
division algebras. The construction we use highlights the
{\it triality} principle, since we put a natural Lie algebra 
structure on the direct sum

\setlength{\unitlength}{4mm}
\begin{picture}(30,10)(-12,2)
\put(-3,6){$\fg (\AA,\BB)=$}
\put(1.4,6){$\AA_1\ot\BB_1$}
\put(4.8,6.2){\line(1,0){2}}
\put(7,6){$\ft(\AA){\times}\ft(\BB)$}
\put(11,7){\line(2,3){1.5}}
\put(12.5,9.5){$\AA_3\ot\BB_3$}
\put(11,5.5){\line(2,-3){1.5}}
\put(12.5,2.5){$\AA_2\ot\BB_2$}
\end{picture}

\noindent where $\ft(\AA)$ is a certain triality algebra associated to $\AA$. 
This structure was actually discovered by Barton and Sudbery
(following suggestions of Ramond), who showed it was equivalent 
to the original construction of Tits \cite{bs}. We give a much more 
direct and simple proof, which was also obtained independently by 
Dadok and Harvey \cite{dh}. 

All this leads to a simple description of the exceptional root systems,
the key point for the dimension formulas being that the roots of 
$\fg (\AA,\OO)$ are naturally partitioned into intervals whose endpoints
are linear functions of $a$. This allows one to explicitly write down
infinite series of formulas generalizing those of Deligne, Cohen and de Man,
see theorem 3.3. For example, specializing to just Cartan powers
of the adjoint representation  we obtain:

\begin{prop}
Let $\fg =\fsl_2, \fsl_3, \fg_2, \fso_8, \ff_4, \fe_6, \fe_7, \fe_8$,
with $a=-4/3,-1,-2/3, 0,1,2,4,8$ respectively. Then
$$\dim \fg^{(k)} = \frac{3a+2k+5}{3a+5} \frac{\binom{k+2a+3}{k}
\binom{k+\frac{5a}{2}+3}{k}\binom{k+3a+4}{k}}
{\binom{k+\frac{a}{2}+1}{k}\binom{k+a+1}{k}}.$$
\end{prop}

The perspective also naturally uncovers the representations distinguished
by Freudenthal and dimension formulas for their
Cartan powers, see theorems 4.3 and 5.3. In particular it leads to new models for the standard 
representations in the second and third rows of Freudenthal's magic square.

In a companion paper to this one \cite{lm4},
we discuss the decomposition formulas of 
Deligne and those of Vogel
(coming from his
conjectured \lq\lq universal Lie algebra\rq\rq , see \cite{vog}) from a geometric perspective. We are able to account
for nearly all the factors that appear in their decompositions using
elementary algebraic geometry. This paper is the fourth in a series
exploring connections between representation theory and the projective
geometry of rational homogeneous varieties (see also \cite{lm1,lm2,lm3}).

\section{Triality and the  Vinberg construction}

For $\AA$ a normed algebra over a field $k$, let 
$$T(\AA)=\{\th=(\th_1,\th_2,\th_3)\in SO(\AA)^3,\;\;
\th_3(xy)=\th_1(x)\th_2(y)\;\forall x,y\in\AA\}.$$
There are three natural actions of
$T(\AA)$ on $\AA$ corresponding to its three projections on $SO(\AA)$, 
and we denote these representations by $\AA_1$, $\AA_2$, $\AA_3$.

If $\AA$ is a real Cayley algebra, it is a classical fact that $T(\AA)$ is 
an algebraic group of type $D_4$. In this case the representations 
$\AA_1$, $\AA_2$, $\AA_3$ are non-equivalent and they are exchanged
by the outer automorphism $t$ of $T(\AA)$ of order 3  defined 
by $t.\th = (\th_2,\th_3,\th_1)$. This is the famous {\sl triality 
principle}, encoded in the triple symmetry of the Dynkin diagram for $D_4$.
For the other   real normed division algebras $\AA$, 
we get the following types for the Lie algebra $\ft(\AA)$  of $T(\AA)$, see 
\cite{bs}:

$$\begin{array}{cccc}
\ft(\RR) & \ft(\CC) & \ft(\HH) & \ft(\OO) \\
 0 & \RR^2 & \fso_3\times\fso_3\times\fso_3 & \fso_8
\end{array}$$

Now let $\AA$ and $\BB$ be two normed algebras. We define on 
$$\fg = \ft(\AA)\times\ft(\BB)\oplus (\AA_1\ot\BB_1)
\oplus (\AA_2\ot\BB_2) \oplus (\AA_3\ot\BB_3)
$$
a $\ZZ_2\times\ZZ_2$-graded Lie algebra structure by the following conditions:
\begin{itemize}
\item $\fg_0=\ft(\AA)\times\ft(\BB)$;
\item the bracket of an element of $\ft(\AA)\times\ft(\BB)$
with one of $\AA_i\ot\BB_i$ is given by the actions of $\ft(\AA)$ 
on $\AA_i$ and $\ft(\BB)$ on $\BB_i$, that is
$$[\th^{\AA},u_i\ot v_i]=\th^{\AA}_i(u_i)\ot v_i, \qquad 
[\th^{\BB},u_i\ot v_i]=u_i\ot \th^{\BB}_i(v_i) ;$$
\item the bracket of two elements in $\AA_i\ot\BB_i$ is given by
the natural map $\La^2(\AA_i\ot\BB_i)=\La^2\AA_i\ot S^2\BB_i\oplus
S^2\AA_i\ot\La^2\BB_i\ra \La^2\AA_i\oplus\La^2\BB_i\ra\ft(\AA)\times
\ft(\BB)$, where the first arrow follows from the quadratic forms
given on $\AA_i$ and $\BB_i$, and the second arrow is dual to the 
map $\ft(\AA)\ra\La^2\AA_i\subset End(\AA_i)$ (and similarly for $\BB)$ 
prescribing the action of $\ft(\AA)$ on $\AA_i$ (which, by
definition, preserves the quadratic form on $\AA_i$). Here duality
is taken with respect to a $\ft(\AA)$-invariant quadratic form on the
reductive algebra $\ft(\AA)$, and the quadratic form on $\La^2\AA_i$
induced by that on $\AA_i$;
\item finally, the bracket of an element of $\AA_i\ot\BB_i$ with one 
of $\AA_j\ot\BB_j$, for $i\ne j$, is given by the following rules,
with obvious notations: 
$$\begin{array}{lcll}
[u_1\ot v_1,u_2\ot v_2]  
 & = & u_1u_2\ot v_1v_2 & \in\AA_3\ot\BB_3, \\
 \hspace{0mm} [u_2\ot v_2,u_3\ot v_3] 
 & = & u_3\Bar{u_2}\ot v_3\Bar{v_2} &
 \in\AA_1\ot\BB_1, \\
 \hspace{0mm} [u_3\ot v_3,u_1\ot v_1] 
  & = & \Bar{u_1}u_3\ot \Bar{v_1}v_3 & 
 \in\AA_2\ot\BB_2.
\end{array}$$
\end{itemize}

\begin{theo} 
This bracket defines a structure of semi-simple Lie algebra on $\fg$, 
whose type is given by Freudenthal's magic square. 
Moreover, each $\fh_i=\ft(\AA)\times\ft(\BB)\oplus\AA_i\ot\BB_i$ is a 
subalgebra of maximal rank of $\fg$.
\end{theo}

Our definition above of the 
Lie bracket on $\fg$ is much simpler  than that in \cite{bs}
since in does not
involve Jordan algebras and their derivations as the Tits construction
does. As a result, below we present  a simpler proof of the fact that $\fg$ is 
indeed a Lie algebra.

The following tables gives the list of possible types for $\fg$ and 
$\fh$. The first table is Freudenthal's magic square. \medskip

\medskip
\begin{center}\begin{tabular}{cccc}
$\fsl_2$ & $\fsl_3$ & $\fsp_6$ & $\ff_4$ \\  
$\fsl_3$ & $\fsl_3\times \fsl_3$ & $\fsl_6$ & $\fe_6$ \\ 
$\fsp_6$ & $\fsl_6$ & $\fso_{12}$ & $\fe_7$  \\
$\ff_4$ & $\fe_6$ & $\fe_7$ & $\fe_8$ 
\end{tabular}\end{center}
\medskip

\medskip
\begin{center}\begin{tabular}{cccclcccc}
 & & $\fsl_2\times\fsl_2\times\fso_8$ & $\fso_9$ \\  
 & & $\CC\times\fsl_2\times\fsl_4$ & $\CC\times\fso_{10}$ \\ 
$\fsl_2\times\fsp_4$ & $\fsl_2\times\CC\times\fsl_4$ &
$\fsl_2\times\fsl_2\times\fso_8$  & $\fsl_2\times\fso_{12}$ \\
$\fso_9$ & $\CC\times\fso_{10}$ & $\fsl_2\times\fso_{12}$ 
& $\fso_{16}$ 
\end{tabular}\end{center}
\medskip

\proof We must check that the Jacobi identity holds in $\fg$.
We begin with a few remarks. 
Denote by $\Psi_i : \La^2\AA_i\ra\ft(\AA)$ the map dual to the 
action of $\ft(\AA)$ on $\AA_i$ with respect to an invariant 
non degenerate quadratic
form $K_{\ft(\AA)}$ on $\ft(\AA)$, and the quadratic form on
$\La^2\AA_i$ induced by the quadratic form $Q=Q_{\AA_i}$ 
on $\AA_i$. We   have
$$K_{\ft(\AA)}(\Psi_i(u\wedge v),\th)=Q(u,\th_i(v))\quad 
\forall u,v\in\AA_i, \forall\th\in\ft(\AA).$$
The action of $\ft(\AA)$ on $\AA_1$ factors through the natural 
representation of $SO(\AA)$, while the actions on $\AA_2$ and $\AA_3$
are induced by the left and right multiplications of $\AA$ on itself. 
More precisely, we have the following formulas:
$$\begin{array}{rcl}
\Psi_1(u\wedge v)_1x & = & Q(u,x)v-Q(v,x)u, \\
\Psi_1(u\wedge v)_2x & = & \Bar{v}(ux)-\Bar{u}(vx), \\
\Psi_1(u\wedge v)_3x & = & (xu)\Bar{v}-(xv)\Bar{u}.
\end{array}$$

(For the case of octonions, these formulas can be deduced from 
\cite{post}, Lecture 15. The other cases are easy.)
Using the compatibility of our construction with the automorphism of 
$\ft(\AA)$ which exchanges the three representations $\AA_i$, we are
reduced to verifying this identity between homogeneous elements in 
the following cases:
 
\begin{enumerate} 
\item $(\ft(\AA),\ft(\AA),\ft(\AA))$-- this is just the Jacobi identity
inside $\ft(\AA)$;

\item $(\ft(\AA),\ft(\AA),\AA_1\ot\BB_1))$-- this case follows from the 
equivariance of the action of $\ft(\AA)$ on $\AA_1$;

\item $(\ft(\AA),\AA_1\ot\BB_1,\AA_1\ot\BB_1))$-- this case follows 
from the equivariance of $\Psi_1$;

\item $(\AA_1\ot\BB_1,\AA_1\ot\BB_1,\AA_1\ot\BB_1))$-- here we must
  check that for $a,b,c,d,e,f\in\AA_1$, 
$$[ [a\ot b, c\ot d], e\ot f] +[ [c\ot d, e\ot f], a\ot b] +
 [ [e\ot f, a\ot b], c\ot d]=0.$$
But the first of these brackets, for example, can be computed as
follows:
$$\begin{array}{rcl}
 \hspace{0mm} [[a\ot b, c\ot d], e\ot f] & = & 
Q(b,d)\Psi_1(a\wedge c)_1e\ot f + Q(a,c)e\ot \Psi_1(b\wedge d)_1f \\
 & = & Q(b,d)Q(a,e)c\ot f - Q(b,d)Q(c,e)a\ot f + \\
 &   & + Q(a,c)Q(b,f)e\ot d - Q(a,c)Q(d,f)e\ot b,
\end{array}$$
and the result easily follows;

\item $(\ft(\AA),\AA_1\ot\BB_1,\AA_2\ot\BB_2))$-- here we need to
  check that 
$$\begin{array}{l}
\hspace{0mm} [ [\th, a\ot b], c\ot d] -[ \th, [a\ot b, c\ot d]] +
 [ a\ot b, [\th, c\ot d]] \\
\hspace*{2cm}
 = [ \th_1(a)\ot b], c\ot d] -[ \th, ac\ot bd] + [ a\ot b, \th_2(c)\ot d]] \\
\hspace*{2cm}
= \left\{ \th_1(a)c-\th_3(ac)+a\th_2(c) \right\}\ot bd =0,
\end{array}$$
and this follows from the infinitesimal triality pronciple for $\th$. 

\item $(\AA_1\ot\BB_1, \AA_1\ot\BB_1, \AA_2\ot\BB_2))$-- here we compute
$$\begin{array}{l}
\hspace{0mm} [ [a\ot b, c\ot d], e\ot f] +[ [c\ot d, e\ot f], a\ot b] +
 [ [e\ot f, a\ot b], c\ot d] \\
\hspace*{2cm}
 = [Q(b,d)\Psi_1(a\wedge c)+Q(a,c)\Psi_1(b\wedge d), e\ot f]
+[ce\ot df, a\ot b]-[ae\ot bf, c\ot d] \\
\hspace*{2cm}
 = Q(b,d)\Psi_1(a\wedge c)_2e\ot f+Q(a,c)e\ot \Psi_1(b\wedge d)_2f
+\Bar{a}(ce)\ot \Bar{b}(df)-\Bar{c}(ae)\ot \Bar{d}(bf).
\end{array}$$
To check that this is zero, we split this expression into its
symmetric and antisymmetric parts with respect to $a$ and $c$. 
To control the symmetric part, we simply let $c=a$, and since 
$\Bar{a}(ae)=Q(a,a)e$, we are left with
$$Q(a,a)e\ot \left\{\Psi_1(b\wedge d)_2f
+\Bar{b}(df)-\Bar{d}(bf)\right\} = 0.$$
Now the antisymmetric part is 
$$ 2Q(b,d)\Psi_1(a\wedge c)_2e\ot f
+\Bar{a}(ce)\ot \Bar{b}(df)-\Bar{c}(ae)\ot \Bar{d}(bf)
-\Bar{c}(ae)\ot \Bar{b}(df)+\Bar{a}(ce)\ot \Bar{d}(bf),$$
which is symmetric in $b$ and $c$. So to check that it vanishes, we
can let $b=d$ and we are left with 
$$ 2Q(b,b)\left\{\Psi_1(a\wedge c)_2e
+\Bar{a}(ce)-\Bar{c}(ae)\right\}\ot f =0.$$

\item $(\AA_1\ot\BB_1, \AA_2\ot\BB_2, \AA_3\ot\BB_3))$-- here we
  compute
$$\begin{array}{l}
\hspace*{0cm}[ [a\ot b, c\ot d], e\ot f] +[ [c\ot d, e\ot f], a\ot b] +
 [ [e\ot f, a\ot b], c\ot d] \\
\hspace*{1cm}
 = [ac\ot bd, e\ot f]+[e\Bar{c}\ot f\Bar{d}, a\ot b]+
[\Bar{a}e\ot \Bar{b}f, c\ot d] \\
\hspace*{1cm}
 = Q(bd,f)\Psi_3(ac\wedge e)+Q(f\Bar{d},b)\Psi_1(e\Bar{c}\wedge a)
+ Q(\Bar{b}f,d)\Psi_2(\Bar{a}e\wedge c)
\end{array}$$
plus a symmetric expression with values in $\ft(\BB)$. But we have, 
$Q(bd,f)=Q(f\Bar{d},b)=Q(\Bar{b}f,d)$, so we just need to check that 
$$\Psi_3(ac\wedge e)+\Psi_1(e\Bar{c}\wedge a)
+ \Psi_2(\Bar{a}e\wedge c)=0.$$
This follows from the triality principle by duality: indeed, for 
every $\th\in\ft(\AA)$, we have
$$\begin{array}{rcl}
K_{\ft(\AA)}(\Psi_3(ac\wedge e),\th) & = & -Q(\th_3(ac),e) \\
 & = & -Q(\th_1(a)c+a\th_2(c),e) \\
 & = & -Q(\th_1(a),e\Bar{c})-Q(\th_2(c),\Bar{a}e) \\
 & = & -K_{\ft(\AA)}(\Psi_1(e\Bar{c}\wedge a),\th)
-K_{\ft(\AA)}(\Psi_2(\Bar{a}e\wedge c),\th),
\end{array}$$
and the result follows. 
\end{enumerate}
This proves that we have endowed $\fg$ with a Lie algebra structure. 
This algebra is reductive. There is a natural quadratic form
${\mathcal Q}$ on $\fg$ defined by the fact that the factors of $\fg$ 
are mutually orthogonal, each one being endowed with its natural 
quadratic form. 

\begin{lemm}
The following nondegenerate quadratic form on $\fg$ is
$\fg$-invariant:
$$K=K_{\ft(\AA)} +K_{\ft(\BB)}+\sum_{i=1}^3Q_{\AA_i}\ot Q_{\BB_i}.$$
\end{lemm}

Since the center of $\fg$ is
trivial, we conclude that $\fg$ is semi-simple. Moreover, any Cartan
subalgebra of $\ft(\AA)\times\ft(\BB)$ will be a Cartan subalgebra 
of $\fg$: in particular, ${\rm rank}(\fg)={\rm rank}(\ft(\AA))+
{\rm rank}(\ft(\BB))$. Finally, knowing the ranks and dimensions of
the semi-simple Lie algebra $\fg$ and its reductive subalgebra
$\fh$, we easily check that 
their types are given by Freudenthal's square and the table below.\qed
\bigskip

The triality Lie algebras can be generalized to $r$-ality for all $r$
to recover the generalized Freudenthal chart (see \cite{lm1}). 
For $r>3$ we have
$$
\ft_r(\RR)=0, \qquad\ft_r(\CC)=\CC^{\op (r-1)}, 
\qquad\ft_r(\HH)=\fsl_2^{\times r}
$$
and
$$
 \fg_r(\AA,\BB) =  \ft_r(\AA) +  \ft_r(\BB) + 
\bigoplus_{1\leq i<j\leq r}\AA_{ij}\ot\BB_{ij}
$$
This model is useful for more generalized dimension formulas, see section 7.

\bigskip

\section{The exceptional series}

From now on we work over the complex numbers.

For $\BB=\OO$, our construction gives the last line of Freudenthal 
square. Let us describe the root system of $\fg$. 
For this we choose Cartan subalgebras of $\fso_8$ and 
$\ft(\AA)$. Their product is a Cartan sublagebra of $\fg$,  and the 
corresponding root spaces in $\fg$ are the root spaces in $\fso_8$ and
$\ft(\AA)$  and   the weight spaces 
of the tensor products $\AA_i\ot\OO_i$. Thus the roots of $\fg$ are
\begin{itemize}
\item the roots of $\fso_8$,
\item the roots of $\ft(\AA)$,
\item the weights $\mu+\nu$, with $\mu$ a weight of $\AA_i$ and $\nu$
a weight of $\OO_i$.
\end{itemize}
To get a set of positive roots we choose  linear forms $l$ and
$l_{\AA}$ on the root lattices, that are strictly positive on positive 
roots. More precisely, we choose $l=l_1\e_1^*+l_2\e_2^*+l_3\e_3^*+l_4\e_4^*$
with $l_1\gg l_2\gg l_3\gg l_4$. (Here and in what follows,
we use the notations and conventions of \cite{bou}.) Then the linear form $ml+l_{\AA}$,
where $m\gg 1$, will be positive on the following set of positive
roots of $\fg$:
\begin{itemize}
\item the positive roots of $\fso_8$,
\item the positive roots of $\ft(\AA)$,
\item the weights $\mu+\nu$, with $\mu$ a weight of $\AA_i$ and $\nu$
a weight of $\OO_i$ such that $l(\nu)>0$. 
\end{itemize}
These weights $\nu$ of $\OO_i$ such that $l(\nu)>0$ are given by the
following tables: \medskip

\begin{center}\begin{tabular}{cp{2cm}cp{2cm}c}
$\OO_1$ & \qquad &  $\OO_2$ & \qquad & $\OO_3$ \\
  & & & & \\
\setlength{\unitlength}{4mm}
\begin{picture}(3,3.3)(0,0)
\put(-.4,1){$1-1$}
\put(1.9,1.5){\line(2,3){.6}}
\put(1.9,1.2){\line(2,-3){.6}}
\put(2.7,2.2){$\frac{1}{2}$}
\put(2.7,.2){$\frac{1}{2}$}
\end{picture}
& & 
\setlength{\unitlength}{4mm}
\begin{picture}(3,3.3)(0,0)
\put(-.4,1){$\frac{1}{2}-1$}
\put(1.9,1.5){\line(2,3){.6}}
\put(1.9,1.2){\line(2,-3){.6}}
\put(2.7,2.2){$\frac{1}{2}$}
\put(2.7,0){$1$}
\end{picture}
& & 
\setlength{\unitlength}{4mm}
\begin{picture}(3,3.3)(0,0)
\put(-.5,1){$\frac{1}{2}-1$}
\put(1.9,1.5){\line(2,3){.6}}
\put(1.9,1.2){\line(2,-3){.6}}
\put(2.7,.2){$\frac{1}{2}$}
\put(2.7,2.2){$1$}
\end{picture}
\\
\setlength{\unitlength}{4mm}
\begin{picture}(3,3.3)(0,0)
\put(-.4,1){$0-1$}
\put(1.9,1.5){\line(2,3){.6}}
\put(1.9,1.2){\line(2,-3){.6}}
\put(2.7,.2){$\frac{1}{2}$}
\put(2.7,2.2){$\frac{1}{2}$}
\end{picture}
& & 
\setlength{\unitlength}{4mm}
\begin{picture}(3,3.3)(0,0)
\put(-.5,1){$\frac{1}{2}-1$}
\put(1.9,1.5){\line(2,3){.6}}
\put(1.9,1.2){\line(2,-3){.6}}
\put(2.7,0){$0$}
\put(2.7,2.2){$\frac{1}{2}$}
\end{picture}
& & 
\setlength{\unitlength}{4mm}
\begin{picture}(3,3.3)(0,0)
\put(-.5,1){$\frac{1}{2}-1$}
\put(1.9,1.5){\line(2,3){.6}}
\put(1.9,1.2){\line(2,-3){.6}}
\put(2.7,2.2){$0$}
\put(2.7,.2){$\frac{1}{2}$}
\end{picture}
\\
\setlength{\unitlength}{4mm}
\begin{picture}(3,3.3)(0,0)
\put(-.4,1){$0-0$}
\put(1.9,1.5){\line(2,3){.6}}
\put(1.9,1.2){\line(2,-3){.6}}
\put(2.7,2.2){$\frac{1}{2}$}
\put(2.7,.2){$\frac{1}{2}$}
\end{picture}
& & 
\setlength{\unitlength}{4mm}
\begin{picture}(3,3.3)(0,0)
\put(-.5,1){$\frac{1}{2}-0$}
\put(1.9,1.5){\line(2,3){.6}}
\put(1.9,1.2){\line(2,-3){.6}}
\put(2.7,2.2){$\frac{1}{2}$}
\put(2.7,0){$0$}
\end{picture}
& & 
\setlength{\unitlength}{4mm}
\begin{picture}(3,3.3)(0,0)
\put(-.5,1){$\frac{1}{2}-0$}
\put(1.9,1.5){\line(2,3){.6}}
\put(1.9,1.2){\line(2,-3){.6}}
\put(2.7,.2){$\frac{1}{2}$}
\put(2.7,2.2){$0$}
\end{picture}
\\
\setlength{\unitlength}{4mm}
\begin{picture}(3,3.3)(0,0)
\put(-.4,1){$0-0$}
\put(1.9,1.5){\line(2,3){.6}}
\put(1.9,1.2){\line(2,-3){.6}}
\put(2.7,.2){$\frac{1}{2}$}
\put(2.7,2.2){$-\frac{1}{2}$}
\end{picture}
& & 
\setlength{\unitlength}{4mm}
\begin{picture}(3,3.3)(0,0)
\put(-.5,1){$\frac{1}{2}-0$}
\put(1.9,1.5){\line(2,3){.6}}
\put(1.9,1.2){\line(2,-3){.6}}
\put(2.7,0){$0$}
\put(2.7,2.2){$-\frac{1}{2}$}
\end{picture}
& & 
\setlength{\unitlength}{4mm}
\begin{picture}(3,3.3)(0,0)
\put(-.5,1){$\frac{1}{2}-0$}
\put(1.9,1.5){\line(2,3){.6}}
\put(1.9,1.2){\line(2,-3){.6}}
\put(2.7,2.2){$0$}
\put(2.7,.2){$-\frac{1}{2}$}
\end{picture}
\end{tabular}\end{center}\medskip

E.g., the first weight in the first 
column is $\a_1+\a_2+\frac 12\a_3+\frac 12\a_4$.

\begin{rem} From this explicit description of the root system of
  $\fg$, it is quite easy to extract the set of simple roots, from
which one can readily obtain the Dynkin diagram of $\fg$. Observe in
particular that if we normalize the invariant scalar product of the 
(dual of the) Cartan algebra of $\fso_8$ in such a way that the root
lengths equal two, then the length of a root of the form $\m+\nu$
equals $(\m,\m)+(\nu,\nu)=1+(\nu,\nu)$. For $\AA\ne\RR$, this is 
larger than one, so it must in fact equal two. This fixes the 
relative normalization of the invariant scalar product on the 
Cartan subalgebra of $\ft(\AA)$, and shows that $\fg$ must be 
simply laced. On the contrary, for $\AA=\RR$ we get roots of length
one ($\fg=\ff_4$ !), and there is no problem of normalization. 
\end{rem}

\begin{prop}

1. With the ordering above, the following are, in order, the three
highest roots of $\fg$.
$$
\setlength{\unitlength}{4mm}
\begin{picture}(7.2,2.5)(-3,.3)
\put(-3,1){$\tilde{\a}_{\OO} =$} 
\put(-.5,1){$1-2$}
\put(1.9,1.5){\line(2,3){.6}}
\put(1.9,1.2){\line(2,-3){.6}}
\put(2.7,2.2){$1$}
\put(2.7,.2){$1$}
\end{picture} \quad 
\setlength{\unitlength}{4mm}
\begin{picture}(7.2,2.5)(-3,.3)
\put(-3,1){$\b_1 =$} 
\put(-.5,1){$1-1$}
\put(1.9,1.5){\line(2,3){.6}}
\put(1.9,1.2){\line(2,-3){.6}}
\put(2.7,2.2){$1$}
\put(2.7,.2){$1$}
\end{picture}, \quad 
\begin{picture}(7.2,2.5)(-3,.3)
\put(-3,1){$\b_2 =$} 
\put(-.5,1){$1-1$}
\put(1.9,1.5){\line(2,3){.6}}
\put(1.9,1.2){\line(2,-3){.6}}
\put(2.7,2.2){$0$}
\put(2.7,.2){$1$}
\end{picture} \quad 
 $$

They are all the simple roots of $\fg$ annhilated by
 the torus of $\ft (\AA)$, in fact the next
highest root is $\b_3=\o_1+\mu^+$ where
$\mu^+$ is the highest weight of $\AA_1$.

2. Any positive weight of $\fg$ annhilated by 
the torus of $\ft(\AA)$ is a linear combination of the following four
weights:
$$
\setlength{\unitlength}{4mm}
\begin{picture}(8.2,2.5)(-5,.3)
\put(-3.7,1){$\o(\fg) =$} 
\put(-.5,1){$1-2$}
\put(1.9,1.5){\line(2,3){.6}}
\put(1.9,1.2){\line(2,-3){.6}}
\put(2.7,2.2){$1$}
\put(2.7,0){$1$}
\end{picture} \quad 
\setlength{\unitlength}{4mm}
\begin{picture}(8.2,2.5)(-5,.3)
\put(-4.5,1){$\o(X_2) =$} 
\put(-.5,1){$2-3$}
\put(1.9,1.5){\line(2,3){.6}}
\put(1.9,1.2){\line(2,-3){.6}}
\put(2.7,2.2){$2$}
\put(2.7,0){$2$}
\end{picture} \quad 
\begin{picture}(8.2,2.5)(-5,.3)
\put(-4.5,1){$\o(X_3) =$} 
\put(-.5,1){$3-4$}
\put(1.9,1.5){\line(2,3){.6}}
\put(1.9,1.2){\line(2,-3){.6}}
\put(2.7,2.2){$2$}
\put(2.7,0){$3$}
\end{picture} \quad 
\begin{picture}(8.2,2.5)(-5,.3)
\put(-4.5,1){$\o(Y_2^*) =$} 
\put(-.5,1){$2-2$}
\put(1.9,1.5){\line(2,3){.6}}
\put(1.9,1.2){\line(2,-3){.6}}
\put(2.7,2.2){$1$}
\put(2.7,0){$1$}
\end{picture} $$

They occur respectively, $\o(\fg)=\o_2$ as the highest weight of
$\fg$, $\o(X_2)=\o_1+\o_3+\o_4$ as the highest weight of $\La^2\fg$, 
$\o(X_3)=2\o_1+2\o_3$ as the highest weight of $\La^3\fg$, 
and $\o(Y_2^*)=2\o_1$ as the highest weight of $S^2\fg - \fg\up 2$. 
 
 3. The half-sum of the positive roots of $\fg$ is $\rho = \rho_{\ft(\AA)}
+\rho_{\ft(\OO)}+a\g_{\ft(\OO)}$, where $\g_{\ft(\OO)}=2\o_1+\o_4$. 
\end{prop}

The names of the weights are borrowed from \cite{cdm}. We
show below the representations are indeed those of \cite{cdm} for
$\fe_8$ below, the other cases are safely left to the reader. 

\proof Everything is clear except for the assertion about $\o (Y_2^*)$, which
is a consequence of the following observations. 

Let $\m$ be 
a weight of $S^2\fg$ such that $2\tilde{\a}\ge\m >2\o_1$. 
Such a weight $\m$ must be the sum of two positive roots $\g$ and
$\d$. Suppose that $\g, \d\ne\tilde{\a}$. Then $\g,\d$ have
coefficients at most one, hence $\m$ has coefficients at most two,
when expressed in terms of simple roots. Since $\m>2\o_1=
2\a_1+2\a_2+\a_3+\a_4$, 
we have   $\m=2\a_1+ 2\a_2+ 2\a_3+ \a_4$ (up to exchanging
$\a_3$ and $\a_4$) hence $\g=\a_1+ \a_2+ \a_3+ \a_4$ and 
$\d=\a_1+ \a_2+ \a_3$. Since in that case $\m-\tilde{\a}$ is not a 
root, this implies that each  possible $\m$ has multiplicity one
inside $S^2\fg$. But it also has multiplicity one inside
the irreducible component of highest weight $2\tilde{\a}$. 

The situation is different for $2\o_1$, whose multiplicity is at least
  $3$ since there are already $3$ different ways to
write it as the sum of two roots of $\fso_8$.
To conclude, we just need to check that the multiplicty of $2\o_1$ is strictly 
larger than its multiplicity inside
the irreducible $\fg$-module of highest weight $2\tilde{\a}$. 
But since $2\o_1$ and $2\tilde{\a}$ both have support on the weight
lattice of $\fso_8$, it follows from Kostant's multiplicity formula 
that this multiplicity can be computed directly in $\fso_8$, where
we check that it is two. We are done. \qed

\begin{rem}
Consider the weights $\o$ of $\fg$ that have support on the Cartan 
subalgebra of $\fso_8$. Obviously, they must belong to the weight 
lattice of $\fso_8$, but there are more conditions imposed by the 
roots of $\fg$ of the form $\m+\nu$, $\nu$ a weight of $\OO_i$, 
$\m$ a weight of $\AA_i$: namely, $2(\o,\nu)/(\m+\nu,\m+\nu)$ must be
an integer. We have $(\m,\m)=1$, and $(\nu,\nu)=1$ as well (except 
in the case where $\AA=\RR$, for which $\nu =0$). Thus our conditions
reduce to $(\o,\nu)\in\ZZ$ for each $\nu$. If we write $\o=o_1\o_1+
o_2\o_2+o_3\o_3+o_4\o_4$, this means that the integers $o_1, o_2, 
o_3, o_4$ must be such that $o_1+o_3, o_1+o_4$ and $o_3+o_4$ are even. 
This defines a sub-lattice of index four of the weight lattice of 
$\fso_8$, and it is straightforward to check that the cone of 
positive weights in this lattice is precisely the cone of non negative
linear combinations of the four weights $\o(\fg)$, $\o(X_2)$, 
$\o(X_3)$ and $\o(Y_2^*)$.
\end{rem} \medskip

\begin{exam} Consider  the case of $e_8$, i.e.,  $\AA=\OO$. 
We denote
the roots and weights
of $\ft(\AA)=\fso_8$
 with primes. We first   determine the set of simple roots of $\fe_8$. 
They must be among the simple roots $\a_i$ of $\ft(\OO)$, 
the simple roots $\a'_j$ of $\ft(\AA)$, and the weights 
$$\g_1=\o_3-\o_4-\o'_1,\quad \g_2=\o_1-\o_4-\o'_3,\quad 
\g_3=\o_1-\o_3-\o'_4,$$
which are the smallest positive roots inside $\AA_1\ot\OO_1, 
\AA_2\ot\OO_2, \AA_3\ot\OO_3$
respectively. But $\a_3=\a_4+2(\o_3-\o_4)=\a_4+2\g_1+2\o'_1$, where 
$2\o'_1$ belongs to the root lattice of $\fso_8$, showing that $\a_3$
cannot be a simple root of $\fg$. Neither can $\a_1$ for the same 
reason. The same conclusion holds for $\g_2$ because of the 
relation $\g_2-\g_1-\g_3=\o'_1+\o'_4-\o'_3=\a'_1+\a'_2+\a'_4$.
Since we know we must have 8 simple roots, they
must be $\a_2,\a_4,\a'_1,\a'_2,\a'_3,
\a'_4,\g_1,\g_3$. Using them, we easily 
deduce  the Dynkin diagram of $e_8$: we have a subdiagram of type $\fso_8$ corresponding to 
$\a'_1,\a'_2,\a'_3,\a'_4$, and we attach to it 4 other nodes according
to the non zero scalar products $(\a_2, \a_4)$, $(\a_4, \g_1)$, 
$(\g_1, \a'_1)$ and $(\g_3, \a'_4)$:

\setlength{\unitlength}{6mm}
\begin{picture}(20,7)(-7,-2.3)
\multiput(0,0)(2,0){5}{$\bullet$}
\multiput(9.5,1.5)(1.5,1.5){2}{$\bullet$}
\multiput(.2,.15)(2,0){4}{\line(1,0){2}}
\multiput(8.2,.2)(1.5,1.5){2}{\line(1,1){1.5}}
\put(9.5,-1.5){$\bullet$}
\put(8.2,0){\line(1,-1){1.5}}
\put(-.2,.5){$\alpha_2$}
\put(1.8,.5){$\alpha_4$}
\put(3.8,.5){$\gamma_1$}
\put(10.3,3.2){$\gamma_3$}
\put(5.8,.6){$\alpha'_1$}
\put(7.6,.6){$\alpha'_2$}
\put(8.7,-1.8){$\alpha'_3$}
\put(8.8,1.8){$\alpha'_4$}
\end{picture}

It is now a simple computation to express the weights $\o(\fg)$, 
$\o(X_2)$, $\o(X_3)$ and $\o(Y_2^*)$ in terms of our simple roots. 
We obtain 

$$
\setlength{\unitlength}{4mm}
\begin{picture}(12,3.5)(1,.3)
\put(-3.7,1){$\o(\fg) =$} 
\put(-.5,1){$2-3-4-5-6$}
\put(3.9,1.5){\line(1,1){.6}}
\put(5.1,2.7){\line(1,1){.6}}
\put(3.9,1.2){\line(1,-1){.6}}
\put(6.9,2.3){$4$}
\put(8,3.6){$2$}
\put(6.9,-.2){$3$}
\put(8,1){$=\o_8(\fe_8),$} 
\end{picture} \qquad
\setlength{\unitlength}{4mm}
\begin{picture}(12,3.5)(-5,.3)
\put(-5,1){$\o(X_2) =$} 
\put(-.5,1){$3-6-8-10-12$}
\put(3.9,1.5){\line(1,1){.6}}
\put(5.1,2.7){\line(1,1){.6}}
\put(3.9,1.2){\line(1,-1){.6}}
\put(7.9,2.3){$8$}
\put(9,3.6){$4$}
\put(7.9,-.2){$6$}
\put(9.5,1){$=\o_7(\fe_8)$} 
\end{picture} $$

$$
\setlength{\unitlength}{4mm}
\begin{picture}(12,3.5)(1,.3)
\put(-4.5,1){$\o(X_3) =$} 
\put(-.5,1){$4-8-10-13-16$}
\put(3.9,1.5){\line(1,1){.6}}
\put(5.1,2.7){\line(1,1){.6}}
\put(3.9,1.2){\line(1,-1){.6}}
\put(8.2,2.3){$11$}
\put(9.6,3.6){$6$}
\put(8.2,-.3){$8$}
\put(9.5,1){$=\o_6(\fe_8),$} 
\end{picture} \qquad
\setlength{\unitlength}{4mm}
\begin{picture}(12,3.5)(-5,.3)
\put(-4.5,1){$\o(Y_2^*) =$} 
\put(-.5,1){$2-4-6-8-10$}
\put(3.9,1.5){\line(1,1){.6}}
\put(5.1,2.7){\line(1,1){.6}}
\put(3.9,1.2){\line(1,-1){.6}}
\put(7.8,2.3){$7$}
\put(8.8,3.6){$4$}
\put(7.8,-.2){$5$}
\put(9.1,1){$=\o_1(\fe_8).$} 
\end{picture} $$

\medskip
This shows our terminology agrees with that of \cite{cdm} in the case
of $\fe_8$.
\end{exam}

\medskip
Now we make a few observations on the weights of $\AA_i$. 
First note that since $\AA_i$ has an invariant quadratic form, 
the set of its weights is symmetric with respect to the origin. 
In particular, their sum is zero. The weight structure is as follows:

\setlength{\unitlength}{4mm}
\begin{picture}(30,11)(-8,-1)
\put(0,0){$\bullet$}
\put(0,1){$\bullet$}
\put(0,2){$\bullet$}
\put(-1,3){$\bullet$}
\put(1,3){$\bullet$}
\put(0,4){$\bullet$}
\put(0,5){$\bullet$}
\put(0,6){$\bullet$}
\put(0.2,0.2){\line(0,1){2}}
\put(0.2,4.2){\line(0,1){2}}
\put(0.2,2.2){\line(1,1){1}}
\put(0.2,2.2){\line(-1,1){1}}
\put(0.2,4.2){\line(1,-1){1}}
\put(0.2,4.2){\line(-1,-1){1}}
\put(0,8){${\mathbb O}$}

\put(7,2){$\bullet$}
\put(6,3){$\bullet$}
\put(8,3){$\bullet$}
\put(7,4){$\bullet$}
\put(7.2,2.2){\line(1,1){1}}
\put(7.2,2.2){\line(-1,1){1}}
\put(7.2,4.2){\line(1,-1){1}}
\put(7.2,4.2){\line(-1,-1){1}}
\put(13,3){$\bullet$}
\put(15,3){$\bullet$}
\put(7,8){${\mathbb H}$}
\put(14,8){${\mathbb C}$}
\put(21,3){$\bullet$}
\put(21,8){${\mathbb R}$}
\end{picture}

In particular, when $\m$ describes the weights of $\AA_i$, the 
integer $(\rho,\mu)$ takes each value in the interval
$[1-\frac{a}{2},\frac{a}{2}-1]$   once (this is the empty 
interval for $a=1$), plus the value zero once more. We call
this set of values $v(\AA)$. 
 
Now look at the inner products of the weights $\o(\fg), \o(X_2), \o(X_3)$ and
$\o(Y_2^*)$ with the positive roots of $\fg$. Since these four weights come 
from $\fso_8$ only, the pairing  is zero on the roots coming from
$\ft(\AA)$. Moreover, on the roots of the form $\m+\nu$, the pairing
 depends only  on $\m$. We get the following possibilities:
\medskip

\begin{center}\begin{tabular}{ccp{2cm}cc}
0122 & (12) & & 1232 & (3$\frac{5}{2}$) \\
1000 & (10) & & 1110 & (2$\frac{1}{2}$) \\
0100 & (10) & & 0110 & (1$\frac{1}{2}$) \\
0120 & (11) & & 0010 & (0$\frac{1}{2}$) \\
1122 & (22) & & 1221 & (3$\frac{3}{2}$) \\
1100 & (20) & & 1121 & (2$\frac{3}{2}$) \\
1120 & (21) & & 0121 & (1$\frac{3}{2}$) \\
1220 & (31) & & 0001 & (0$\frac{1}{2}$) \\
1242 & (33) & & 1231 & (32) \\
1222 & (32) & & 1111 & (21) \\
1342 & (43) & & 0111 & (11) \\
2342 & (53) & & 0011 & (01) 
\end{tabular}\end{center}

\medskip
The first column comes from the positive roots of $\fso_8$, each
possibility occurs 
exactly once. The second column comes from the weights $\mu$ of the 
three modules $\OO_i$: we denote by $\Sigma$ the set of these
weights. Here each possibility occurs for exactly $a$ 
positive roots of $\fg$. In parenthesis, we have also given the values 
$u$ of $(\rho_{\ft(\OO)},\a)$ and $v$ of $(\g_{\ft(\OO)},\a)$. For the first
column, this means that  $(\rho,\a)=u+av$. For the 
second column, the values taken by $\rho$ on the $a$ positive 
roots for each case will be the set $v(\AA)$ translated by $u+av$. 
This is the information we need to apply the Weyl dimension formula. 

\begin{theo}\label{super}
 The dimension of the irreducible $\fg$-module
with highest weight $\o= p\o(\fg)+q\o(X_2)+r\o(X_3)+s\o(Y_2^*)$
is given by the following formula:
$${\rm dim}\;V_{\o} = \prod_{\a\in\Delta_+(\fso_8)\cup\Sigma}
\frac{(a\g_{\ft(\OO)}+\rho_{\ft(\OO)}+\o,\a)}
{(a\g_{\ft(\OO)}+\rho_{\ft(\OO)},\a)}
\prod_{\beta\in\Sigma}
\frac{\begin{pmatrix}
(a\g_{\ft(\OO)}+\rho_{\ft(\OO)}+\o,\b)+\frac{a}{2}-1 \\ (\o,\b)
\end{pmatrix}}
{\begin{pmatrix}
(a\g_{\ft(\OO)}+\rho_{\ft(\OO)}+\o,\b)-\frac{a}{2} \\ (\o,\b)
\end{pmatrix}} .
$$

For each choice of $p,q,r,s$, this formula gives a rational 
function of $a$, whose numerator and denominator are products 
of $6p+12q+16r+10s+24$ linear forms.
\end{theo}

\medskip
This formula includes and provides a wide generalization of 15 of the
25 dimension formulas of \cite{cdm}. Since it applies to actual nontrivial
irreducible
representations of $\fd_4,\ff_4,\fe_6,\fe_7,\fe_8$, one could not hope
to apply it to their representations that are zero, negative or
reducible (i.e., two copies of the same representation) for one
of these algebras. When one removes such representations from the
list of 25, only the 15 we are able to account for remain, so in
that sense this is the best possible formula. 

The formula can be made more explicit as follows. 
Each term $abcd \;(uv)$ in the
table above contributes to the product a term $(x+u+av)/(u+av)$,
where $x=ap+bq+cr+ds$. If it is the term is from the second column, 
it also contributes 
$$\frac{(u+av+\frac{a}{2})\cdots (u+av+\frac{a}{2}+x-1)}
{(u+av+1-\frac{a}{2})\cdots (u+av-\frac{a}{2}+x)} =
\frac{(u+av+x-\frac{a}{2}+1)\cdots (u+av+x+\frac{a}{2}-1)}
{(u+av+1-\frac{a}{2})\cdots (u+av+\frac{a}{2}-1)},$$
where the numerator and denominator of the rational function on the
left are products of $x$ linear forms, and those of the rational
function on the right are products of $a-1$ linear forms. 

Specializing to multiples of the highest root, we obtain the 
formula of proposition 1.1 of the introduction. We have so far 
proved this proposition only for $a\ge 0$, but we give a second proof
in section 6  that is valid for the entire series. 

In    Deligne's notations, $\fg^{(k)}$ is $Y_k$. 
Using his parameter $\l$ we get
$$\dim Y_k = \frac{(2k-1)\l-6}{k!\l^k(\l+6)} \prod_{j=1}^k
\frac{((j-1)\l-4)((j-2)\l-5)((j-2)\l-6)}{(j\l-1)((j-1)\l-2)}.$$

Note that the $q$-analogs of our formulas (see e.g.
\cite{ram}, p. 102) are immediate consequences of our methods. 
For example, 
$${\rm dim}_q \fg^{(k)} = \frac{1-q^{3a+2k+5}}{1-q^{3a+5}} 
\frac{{\begin{bmatrix} k+2a+3 \\ k \end{bmatrix}}_q
{\begin{bmatrix} k+\frac{5a}{2}+3 \\ k \end{bmatrix}}_q
{\begin{bmatrix} k+3a+4 \\ k \end{bmatrix}}_q}
{{\begin{bmatrix} k+\frac{a}{2}+1 \\ k \end{bmatrix}}_q
{\begin{bmatrix} k+a+1 \\ k \end{bmatrix}}_q},$$
where ${\begin{bmatrix} k+l \\ k \end{bmatrix}}_q=\frac{(1-q^{l+1})
\cdots (1-q^{l+k})}{(1-q)\cdots (1-q^{k})}$ is the usual Gauss
polynomial.\bigskip

The closed $G$-orbit inside $\PP\fg$, which we call the adjoint
variety and denote by $X_{ad}$, has dimension $6a+9$. From the 
above proposition we can deduce a uniform formula for its degree:

$${\rm deg}\, X_{ad} =
2\frac{(\frac{a}{2}+1)!(a+1)!(6a+9)!}{(\frac{5a}{2}+3)!(2a+3)!(3a+5)!}.
$$

After Freudenthal \cite{fr}, we respectively
 call $X_{F-planes}$, $X_{F-lines}$ and $X_{F-points}$
the closed orbits in $\PP X_2$, $\PP X_3$ and $\PP Y_2^*$.
Specializing to the Cartan powers, theorem \ref{super} 
gives the Hilbert functions of these varieties, respectively

$$\begin{array}{rcl}
\dim X_2^{(k)} & = & \frac{k+a+1}{a+1}\frac{k+a+2}{a+2}
\frac{2k+2a+3}{2a+3}\frac{3k+3a+4}{3a+4}\frac{3k+3a+5}{3a+5}\times 
\\ & & \hspace{3cm} \times
\frac{\binom{k+\frac{3a}{2}+1}{k}\binom{k+\frac{3a}{2}+2}{k}
\binom{k+2a+1}{k}\binom{k+2a+2}{k}\binom{2k+\frac{5a}{2}+3}{2k}
\binom{2k+3a+3}{2k}}
{\binom{k+1}{k}\binom{k+\frac{a}{2}}{k}\binom{k+\frac{a}{2}+1}{k}
\binom{2k+a+2}{2k}\binom{2k+\frac{3a}{2}+2}{2k}}, \\
 & & \\
\dim X_3^{(k)} & = & \frac{2k+\frac{3a}{2}+1}{\frac{3a}{2}+1}
\frac{2k+\frac{3a}{2}+2}{\frac{3a}{2}+2}
\frac{2k+\frac{3a}{2}+3}{\frac{3a}{2}+3}
\frac{4k+3a+3}{3a+3}\frac{4k+3a+4}{3a+4}\frac{4k+3a+5}{3a+5}\times \\
 & & 
\hspace{-1cm} \times
\frac{\binom{k+a}{k}\binom{k+a+1}{k}\binom{k+a+2}{k}
\binom{k+\frac{3a}{2}-1}{k}\binom{k+\frac{3a}{2}}{k}
\binom{k+\frac{3a}{2}+1}{k}
\binom{2k+2a+1}{2k}\binom{2k+2a+2}{2k}\binom{2k+2a+3}{2k}
\binom{3k+\frac{5a}{2}+3}{3k}\binom{3k+3a+2}{3k}}
{\binom{k+1}{k}\binom{k+2}{k}\binom{k+\frac{a}{2}-1}{k}
\binom{k+\frac{a}{2}}{k}\binom{k+\frac{a}{2}+1}{k}
\binom{2k+a}{2k}\binom{2k+a+1}{2k}\binom{2k+a+2}{2k}
\binom{3k+\frac{3a}{2}+3}{3k}\binom{3k+2a+2}{3k}}, \\
 & & \\
\dim Y_2^{*(k)} & = & \frac{2k+\frac{5a}{2}+3}{\frac{5a}{2}+3}
\frac{\binom{k+2a}{k}\binom{k+2a+1}{k}\binom{k+2a+3}{k}
\binom{k+\frac{5a}{2}+2}{k}\binom{2k+3a+5}{2k}}
{\binom{k+\frac{a}{2}-1}{k}\binom{k+\frac{a}{2}+1}{k}
\binom{k+\frac{a}{2}+2}{k}\binom{k+a+1}{k}\binom{k+a+3}{k}
\binom{2k+2a}{2k}},
\end{array}$$
one recovers that $\dim X_{F-planes}= 9a+11$, $\dim X_{F-lines}=
11a+9$, $\dim X_{F-points}=9a+6$, and that their degrees are
$$\deg X_{F-planes}  =  
2^{3a+3}3^2\frac{(9a+11)! a! \frac{a}{2}!(\frac{a}{2}+1)!}
 {(\frac{3a}2+1)!(2a+1)!(2a+3)!(\frac{5a}2+3)!(3a+5)!},$$

$$\deg X_{F-lines} = 
2^{3a+6}3^{2a}\frac{(11a+9)!(\frac{a}{2}-1)!
\frac{a}{2}!(\frac{a}{2}+1)!}
{(\frac{3a}{2}-1)!(\frac{3a}{2}+1)!(2a+3)!(\frac{5a}2+3)!}, $$

$$\deg X_{F-points} = 
2^{a+6}\frac{(6a+9)!(\frac{a}{2}+1)!(\frac{a}{2}+2)!
(a+1)!(a+3)!}{(2a+1)!(2a+3)!(\frac{5a}2+2)!(3a+5)!}.$$

\bigskip
\section{The sub-exceptional series}

In this section we let $\BB=\HH$. Then $\ft(\BB)\simeq\fso_3\times
\fso_3\times\fso_3\simeq\fsl_2\times\fsl_2\times\fsl_2$. 
Note that $\ft(\BB)$ can be naturally identified with 
$Im(\HH)^{\op 3}$, acting on $\HH_1\times\HH_2\times\HH_3$ by 
$$(a,b,c)\mapsto (L_b-R_c, L_c-R_a, L_a-R_b),$$
where $L_a, R_a$ denote the operators of left and right
multiplication by $a$, respectively (see \cite{bs}).
This means that if we denote by $U_1, U_2, U_3$ the
natural 2-dimensional representations of our three copies of 
$\fsl_2$, then $\HH_1=U_2\ot U_3$, $\HH_2=U_3\ot U_1$ and 
$\HH_3=U_1\ot U_2$. Therefore, the roots of $\fg$ are 
\begin{itemize}
\item the roots $\pm\a^1$, $\pm\a^2$, $\pm\a^3$ of 
$\fsl_2\times\fsl_2\times\fsl_2$,
\item the roots of $\ft(\AA)$, 
\item the weights $\pm\frac{1}{2}\a^i\pm\frac{1}{2}\a^j+\m$, 
where $\m$ is a weight of $\AA_k$ and $\{i,j,k\}=\{1,2,3\}$.
\end{itemize}
To get a set of positive roots we choose  linear forms $l$ and
$l_{\AA}$ on the root lattices, that are strictly positive on positive 
roots. More precisely, we choose $l=l_1\a^{1*}+l_2\a^{2*}+l_3\a^{3*}$
with $l_1\gg l_2\gg l_3$. Then the linear form $ml+l_{\AA}$,
where $m\gg 1$, will be positive on the following set of positive
roots of $\fg$:
\begin{itemize}
\item $\a^1$, $\a^2$, $\a^3$,
\item the positive roots of $\ft(\AA)$,
\item the weights $\frac{1}{2}\a^i\pm\frac{1}{2}\a^j+\m$, 
where $\m$ is a weight of $\AA_k$, with $\{i,j,k\}=\{1,2,3\}$ and $i<j$.
\end{itemize}
 
An important difference with the exceptional series is that we 
have a nice geometric model for one of the distinguished
$\fg$-modules $V=\AA_1\ot U_1\op\AA_2\ot U_2\op\AA_3\ot U_3\op
U_1\ot U_2\ot U_3$ .

\begin{theo} There is a natural structure of $\fg$-module 
on $$\setlength{\unitlength}{4mm}
\begin{picture}(20,8.5)(-2,3)
\put(-1,6){$V\; =$}
\put(2.1,6){$\AA_1\ot U_1$}
\put(5.6,6.2){\line(1,0){2}}
\put(8.2,6){$U_1\ot$}
\put(9.9,5){$U_2$}
\put(9.9,7){$U_3$}
\put(11,7.7){\line(2,3){1.5}}
\put(12.5,10.2){$\AA_3\ot U_3$}
\put(11,4.8){\line(2,-3){1.5}}
\put(12.5,1.8){$\AA_2\ot U_2$}
\end{picture}$$
This $\fg$-module $V$ is simple of dimension $6a+8$.\end{theo}

\proof We define the action of $\fg$ on $V$ as follows.
There is already a natural action of the subalgebra $\ft(\AA)\times
\ft(\BB)$, and up to the ternary symmetry we just need to define an
action of $\HH_1\ot\AA_1=U_2\ot U_3\ot\AA_1$. This action is provided 
by the natural maps 
$$\begin{array}{lcl} (U_2\ot U_3\ot\AA_1)\ot (U_1\ot U_2\ot U_3)
& \ra & U_1\ot \AA_1, \\
(U_2\ot U_3\ot\AA_1)\ot (U_1\ot \AA_1)
& \ra & U_1\ot U_2\ot U_3, \\
(U_2\ot U_3\ot\AA_1)\ot (U_2\ot \AA_2)
& \ra & U_3\ot \AA_3, \\
(U_2\ot U_3\ot\AA_1)\ot (U_3\ot \AA_3)
& \ra & U_2\ot \AA_2,
\end{array}$$
which are easily defined using the invariant quadratic forms on 
$U_1,U_2,U_3$ and $\AA_1$, and the natural multiplication map
$\AA_1\ot\AA_2\ra\AA_3$. The verification that this defines a
module structure over the algebra $\fg$ is just a 
computation. \qed

\medskip
The natural $\fg$-invariant symplectic form
$\Omega$ on $V$ may be written
$$\Omega = \Omega_{U_1\ot U_2\ot U_3} + \sum_{i=1}^3\Omega_i,$$
where $\Omega_{U_1\ot U_2\ot U_3}$ is just the tensor product of the
determinants on $U_1$, $U_2$, $U_3$, and $\Omega_i$ is the symplectic
form on $U_i\ot\AA_i$ induced by the determinant on $U_i$ and the
quadratic form $Q_i$ on $\AA_i$.

\begin{prop}

1. With the ordering above,    the three
highest roots of $\fg$ are $\o(\fg)=\tilde{\a}=\a^1=2\o^1$,
  $\o(V)=\o^1+\o^2+\o^3$ and $\o^1+\o^2-\o^3$.

They are all the simple roots of 
$\fg$ annhilated by the torus of $\ft (\AA)$, in fact the next
highest root is $\o^1 +\mu^+$ where
$\mu^+$ is the highest weight of $\AA_1$.

2. Any positive weight of $\fg$ annhilated by 
the torus of $\ft(\AA)$ is a linear combination of the following three
weights: $\o(\fg)=\tilde{\a}=\a^1=2\o^1$,
 $\o(V)=\o^1+\o^2+\o^3$ and $\o(V_2)=2\o^1+2\o^2$

They occur respectively as
the highest weight of $\fg$,   $V$, 
 and $\La^2V$.
 
 3. The half-sum of the positive roots of $\fg$ is $\rho = \rho_{\ft(\AA)}
+\rho_{\ft(\HH)}+a\g_{\ft(\HH)}$, where $\g_{\ft(\HH)}=2\o^1+\o^2$.
\end{prop}

The values of the pairings of
 the weights $\o(\fg), \o(V)$ and  $\o(V_2)$ with the 
positive roots of $\fg$ are obtained as follows. Since these three 
weights come from $\fsl_2\times\fsl_2\times\fsl_2$,
their value is zero on the roots coming from $\ft(\AA)$. Moreover, 
on the roots of the form $\frac{1}{2}\a^i\pm\frac{1}{2}\a^j+\m$, 
their values do not depend on $\m$. We get the following possibilities:

$$\begin{array}{rccrc}
212 & (12) & \hspace*{3cm} & 112 & (1\frac{3}{2}) \\
012 & (11) & & 100 & (0\frac{1}{2}) \\
010 & (10) & & 111 & (11) \\
     &     & & 101 & (01) \\
     &     & & 011 & (1\frac{1}{2}) \\ 
     &     & & 001 & (0\frac{1}{2}) 
\end{array}$$

\bigskip
The first column comes from the positive roots of $\fsl_2\times
\fsl_2\times\fsl_2$, each possibility occurs exactly once. 
The second column comes from the weights of the 
three modules $\HH_i$: we denote theis set by $\Gamma$. 
Each possibility occurs for exactly $a$ 
positive roots of $\fg$. In parenthesis are   the values 
of the parings with $\rho_{\ft(\HH)}$ and $\g_{\ft(\HH)}$. Applying the Weyl 
dimension formula as above we get the following result. 

\begin{theo} The dimension of the irreducible $\fg$-module
with highest weight $\o = p\o(\fg)+q\o(V)+r\o(V_2)$
is given by the following function:
$${\rm dim}\;V_{\o} = \prod_{\a\in\Delta_+(\fsl_2^3)\cup\Gamma}
\frac{(a\g_{\ft(\HH)}+\rho_{\ft(\HH)}+\o,\a)}
{(a\g_{\ft(\HH)}+\rho_{\ft(\HH)},\a)}
\prod_{\beta\in\Gamma}
\frac{\begin{pmatrix}
(a\g_{\ft(\HH)}+\rho_{\ft(\HH)}+\o,\b)+\frac{a}{2}-1 \\ (\o,\b)
\end{pmatrix}}
{\begin{pmatrix}
(a\g_{\ft(\HH)}+\rho_{\ft(\HH)}+\o,\b)-\frac{a}{2} \\ (\o,\b)
\end{pmatrix}} .
$$

For each choice of $p,q,r$, this formula gives a rational 
function of $a$, whose numerator and denominator are products 
of $4p+3q+6r+9$ linear forms.
\end{theo}

\begin{coro}
Let $V$ be the distinguished module, of dimension $6a+8$, of 
a semi-simple Lie algebra $\fg$ in the subexceptional series, 
with $a=-\frac{2}{3}, 0, 1, 2, 4, 8$. Then
$$\begin{array}{rcl}
\dim \fg^{(k)} & = & \frac{2k+2a+1}{2a+1}\frac{
\binom{k+\frac{3a}{2}-1}{\frac{3a}{2}-1}
\binom{k+\frac{3a}{2}+1}{\frac{3a}{2}+1}\binom{k+2a}{2a}}{
\binom{k+\frac{a}{2}-1}{\frac{a}{2}-1} 
\binom{k+\frac{a}{2}+1}{\frac{a}{2}+1}}, \\
\dim V^{(k)} & = & \frac{2a+2k+2}{a+1}\frac{\binom{k+2a+1}{2a+1}
\binom{k+\frac{3a}{2}+1}{\frac{3a}{2}+1}}{\binom{k+
\frac{a}{2}+1}{\frac{a}{2}+1}}, \\
\dim  V_2^{(k)} & = & \frac{ ( 4k+3a+2)}{( k+1)(3a+2)}
\frac{\binom{ k+a+1}{a+1}\binom{ k+\frac{3a}{2}}{\frac{3a}{2}}
\binom{ 2k+2a+1}{2a+1}
\binom{k+\frac{3a}{2}-1}{\frac{3a}{2}-1}
\binom{k+a}{a}}
{\binom{ k+\frac{a}{2}}{\frac{a}{2}}
\binom{ 2k+a}{a}\binom{k+\frac{a}{2}-1}{\frac{a}{2}-1}}.
\end{array}$$
\end{coro}

Let $X\subset\PP V$, $X_{ad}\subset\PP\fg$, $X_{F-planes}\subset\PP V_2$
denote the closed orbits. We recover from the Hilbert functions above that
$\dim X_{ad}= 4a+1$, $\dim X= 3a+3$,
$\dim X_{F-lines}= 5a+2$ and

$$
\deg X_{ad}= \frac{(4a+1)!2(\frac a2 -1)!(\frac a2 + 1) !}
{(2a+1)(\frac{3a}2-1)!(\frac{3a}2+1)!(2a)!}
$$

$$
\deg X=\frac{2 (3a+3)!(\frac a2 +1)!}{ (2a+1)!(\frac{3a}2+1)!}
$$

$$
\deg X_{F-lines}= \frac{(5a+2)!2^{a+3}(\frac a2)!(\frac a2 -1)!}
{(3a+2)!(a+1)!(2a+1)!(\frac {3a}2)!(\frac{3a}2-1)!}.
$$
 
\bigskip
\section{The Severi series}

Now we let $\BB=\CC\ot_{\RR}  \CC=\CC\op\CC$, which
is naturally   the plane $\Pi\subset\CC^3$ of 
equation $z_1+z_2+z_3=0$, acting diagonally on $\CC_1\times\CC_2
\times\CC_3$.  There is a natural identification of $\BB_i$ with 
$\CC_j\ot \CC_k^{-1}\op \CC_j^{-1}\ot \CC_k$ for
$\{i,j,k\}=\{1,2,3\}$, as $\ft(\BB)$-modules. 

Let us denote by $\o_1,\o_2,\o_3$ the highest weights of the 
action of $\ft(\BB)$ on $\CC_1$, $\CC_2$, $\CC_3$, which are
subject to the relation $\o_1+\o_2+\o_3=0$. Then the roots
of $\fg$ are:
\begin{itemize}
\item the roots of $\ft(\AA)$, 
\item the weights $\pm(\o_j-\o_k) +\m$,
where $\m$ is a weight of $\AA_i$ and $\{i,j,k\}=\{1,2,3\}$.
\end{itemize}
To get a set of positive roots we choose  linear forms $l$ and
$l_{\AA}$ on the root lattices, that are strictly positive on positive 
roots. More precisely, we choose $l=l_1\o_{1*}+l_2\o_{2*}$
with $l_1\gg l_2\gg 0$. Then the linear form $l+l_{\AA}$,
will be positive on the following set of positive roots of $\fg$:
\begin{itemize}
\item the positive roots of $\ft(\AA)$,
\item the weights $\o_j-\o_k +\m$,
where $\m$ is a weight of $\AA_i$, $j<k$ and $\{i,j,k\}=\{1,2,3\}$.
\end{itemize}

\smallskip
As for the sub-exceptional series, we have a nice geometric model
for the distinguished 
$\fg$-modules:

\begin{theo} There is a natural structure of $\fg$-module 
on $$\setlength{\unitlength}{4mm}
\begin{picture}(20,8)(-2,3)
\put(-1,6){$W\; =$}
\put(2,6){$\AA_1\ot\CC_1^{-1}$}
\put(5.8,6.2){\line(1,0){2}}
\put(8.2,6){$\CC_1^2\op$}
\put(9.9,5){$\CC_2^2$}
\put(9.9,7){$\CC_3^2$}
\put(10.9,7.9){\line(2,3){1}}
\put(12,9.7){$\AA_3\ot\CC_3^{-1}$}
\put(11,4.7){\line(2,-3){1}}
\put(12,2.4){$\AA_2\ot\CC_2^{-1}$}
\end{picture}$$
This $\fg$-module $W$ is simple of dimension 
$3a+3$.\end{theo}

\proof We just need to define the action of a typical factor 
$\AA_1\ot\CC_2\ot\CC_3^{-1}$ on $W$. This action is given by the 
natural maps
$$\begin{array}{lcl} 
(\AA_1\ot\CC_2\ot\CC_3^{-1})\ot \CC_3^2
& \ra & \AA_1\ot\CC_2\ot\CC_3=\AA_1\ot\CC_1^{-1}, \\
(\AA_1\ot\CC_2\ot\CC_3^{-1})\ot (\AA_1\ot\CC_1^{-1})
& \ra & \CC_1^{-1}\ot\CC_2\ot\CC_3^{-1}=\CC_2^2, \\
(\AA_1\ot\CC_2\ot\CC_3^{-1})\ot (\AA_2\ot\CC_2^{-1})
& \ra & \AA_3\ot\CC_3^{-1},
\end{array}$$
where we use for the first two arrows the fact that
$\CC_1\ot\CC_2\ot\CC_3$ is a trivial $\ft(\BB)$-module, 
for the second arrow the quadratic form on $\AA_1$, and for the 
last arrow the multiplication map $\AA_1\ot\AA_2\ra\AA_3$. 
The action on the other factors is equal to zero. 
We leave to the reader the computations that are  necessary 
to check that this is indeed a Lie algebra action of $\fg$. 
The fact that we get a simple module is obvious. \qed

\medskip Note that there is no natural  symplectic or quadratic form 
on $W$, but a very simple $\fg$-invariant cubic form given by 
$$C(x_1,x_2,x_3,X_1,X_2,X_3) = x_1x_2x_3 + \theta (X_1,X_2,X_3).$$

Here $\theta :\AA_1\ot\AA_2\ot\AA_3\ra\CC$ 
is the triality map, see \cite{bs}.

\begin{prop}
The highest root of $\fg$ is $\o(\fg)=\tilde{\a}=\o_1-\o_3+\m_2$, 
where $\m_2$ is the highest weight of $\AA_2$.

The highest weight of $W$ is $\o(W)=2\o_1$, its lowest weight is 
$-\o(W^*)=2\o_3$. 

The half-sum of the positive roots of $\fg$ is $\rho = \rho_{\ft(\AA)}
+a\g_{\ft(\CC)}$, where $\g_{\ft(\CC)}=\o_1-\o_3$.
\end{prop}

\begin{exam} Let us treat in detail the case where $\AA=\OO$, 
leading to $e_6$ and its minimal representation. The simple 
roots of $\fg$ are those of $\ft(\OO)=\fso_8$, say $\a'_1$, 
$\a'_2$, $\a'_3$, $\a'_4$, and $\g_1=\o_1-\o_2-\m_3$, 
$\g_2=\o_2-\o_3-\m_1$, where $\m_1=\o'_1$, $\m_3=\o'_4$ denote 
the highest weights of $\OO_1$, $\OO_3$ respectively. We get the
following Dynkin diagram:

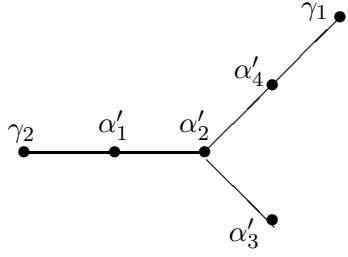
\setlength{\unitlength}{6mm}
\begin{picture}(20,6)(-5,-2.3)
\multiput(4,0)(2,0){3}{$\bullet$}
\multiput(9.5,1.5)(1.5,1.5){2}{$\bullet$}
\multiput(4.2,.15)(2,0){2}{\line(1,0){2}}
\multiput(8.2,.2)(1.5,1.5){2}{\line(1,1){1.5}}
\put(9.5,-1.5){$\bullet$}
\put(8.2,0){\line(1,-1){1.5}}
\put(3.8,.5){$\gamma_2$}
\put(10.3,3.2){$\gamma_1$}
\put(5.8,.6){$\alpha'_1$}
\put(7.6,.6){$\alpha'_2$}
\put(8.7,-1.8){$\alpha'_3$}
\put(8.8,1.8){$\alpha'_4$}
\end{picture}

It is then straightforward to compute $\o(W)$ and $\o(W^*)$ in 
terms of the simple roots. We obtain

$$
\setlength{\unitlength}{4mm}
\begin{picture}(12,3)(1,.3)
\put(-.5,1){$\o(W) =$} 
\put(3.5,1){$\frac{2}{3}-\frac{4}{3}-2$}
\put(3.9,1.5){\line(1,1){.6}}
\put(5.1,2.7){\line(1,1){.6}}
\put(3.9,1.2){\line(1,-1){.6}}
\put(8.2,2.3){$\frac{5}{3}$}
\put(9.3,3.6){$\frac{4}{3}$}
\put(8.2,-.2){$1$}
\put(9.5,1){$=\o_1(\fe_6),$} 
\end{picture} 
\hspace{2cm}
\setlength{\unitlength}{4mm}
\begin{picture}(12,3)(1,.3)
\put(-.5,1){$\o(W^*) =$} 
\put(3.5,1){$\frac{4}{3}-\frac{5}{3}-2$}
\put(3.9,1.5){\line(1,1){.6}}
\put(5.1,2.7){\line(1,1){.6}}
\put(3.9,1.2){\line(1,-1){.6}}
\put(8.2,2.3){$\frac{4}{3}$}
\put(9.3,3.6){$\frac{2}{3}$}
\put(8.2,-.2){$1$}
\put(9.5,1){$=\o_6(\fe_6).$} 
\end{picture} $$

\end{exam}

\medskip
Since the highest root of $\fg$ does depend on $\AA$, we will not 
obtain any rational expression in $a$ for the dimension of $\fg$ 
and its Cartan powers using this model. However, we do obtain such a formula for 
the irreducible $\fg$-modules whose highest weights are linear
combinations of $\o(W)$ and $\o(W^*)$. 

For this we need to compute the values of $\o(W)$ and $\o(W^*)$
on the positive roots of $\fg$. These values are   zero on 
the roots coming from $\ft(\AA)$. 
To compute the other ones, we consider on $\Pi$ the restriction 
of the canonical metric on $\CC^3$. Computing the dual metric we get 
$(\o_i,\o_i)=1/3$ and $(\o_i,\o_j)=-1/6$ for $1\le i\ne j\le 3$. 
It is then straightforward to apply Weyl's dimension formula
and obtain:

\begin{theo} The dimension of the irreducible $\fg$-module
with highest weight $\o = p\o(W)+p^*\o(W^*)$
is given by the following function:
$$\begin{array}{rcl}
{\rm dim}\;V_{\o} & = & \prod_{i\in v(\AA)}
\frac{p+\frac{a}{2}+i}{\frac{a}{2}+i}
\frac{p+p^*+a+i}{a+i}\frac{p^*+\frac{a}{2}+i}{\frac{a}{2}+i} \\
 & = & \frac{(2p+a)(p+p^*+a)(2p^*+a)}{a^3}
\frac{\binom{p+a-1}{p}\binom{p+p^*+\frac{3a}{2}-1}{p+p^*}
\binom{p^*+a-1}{p^*}}{\binom{p+p^*+\frac{a}{2}}{p+p^*}}.
\end{array}$$
\end{theo}

\bigskip

\section{Other models for the exceptional series}

There exist two other  uniform models for the exceptional Lie 
algebras similar to the constructions we used in section 2. In
section 2 we exploited on the triality phenomenon, which is reflected in the 
threefold symmetry of the Dynkin diagram of $\fso_8$. For our 
two other series, we use the simplest Dynkin diagram with twofold
symmetry, which is that of $\fsl_3$, and the simplest one with 
``onefold symmetry'', which is that of $\fsl_2$. This leads to the
three series

$$\begin{array}{rccccccc}
\fg (\AA) & =  & \fso_8 & \op & \ft(\AA) & \op & \OO_1\ot\AA_1
\; \op \;\OO_2\ot\AA_2\;\op \;\OO_3\ot\AA_3  \\ 
\fg (\AA) & = & \fsl_3 & \op & \fsl_3(\AA) & \op & \CC^3\ot\fh_3(\AA)
\;\op \; \CC^{3*}\ot\fh_3(\AA)^*  \\ 
\fg (\AA) & = & \fsl_2 & \op &\fsp_6(\AA) & \op & \CC^2\ot\La^{\langle
  3\rangle}\AA^6
\end{array}$$

Here $SL_3(\AA)$ respectively denotes the Lie groups $Id, ?, {\mathfrak S}_3$
and the four groups on the second row of Freudenthal's magic chart.
$\fh_3(\AA)$  denotes the Jordan algebra over $\AA$ in the last four
cases and $\emptyset$, homotheties, and diagonal $3\times 3$ matrices
in the first three cases (see \cite{lm1}).  Similarly,
$\fsp_6(\AA)$ repspectively denotes $0, \fsl_2, \fsl_2^{\op 3}$ and
the Lie algebras appearing in the third row of Freudenthal's chart.
$\La^{\langle
  3\rangle}\AA^6$ respectively denotes $0$, $S^3\CC^2$ and the subexceptional
  representations $V$.

\medskip
These series show the same remarkable uniformity properties in the 
distributions of the root heights necessary for  nice dimension formulas. 
 The  formulas one obtains only  concern
representations whose highest weights are supported on the weight
lattice of the fixed subalgebra of each series, namely $\fso_8$, 
$\fsl_3$ and $\fsl_2$ respectively. The rank of this subalgebra 
is maximal for the first series so we won't be able to extract more
information from the other two series. 

\smallskip
Let's consider, nevertheless, our models in the second series, involving the 
action of $\fsl_3(\AA)$ on the Jordan algebra $\fh_3(\AA)$. 
A natural Cartan subalgebra of $\fg(\AA)$ is obtained as the 
direct sum of Cartan subalgebras of $\fsl_3$ and $\fsl_3(\AA)$. 
We choose a linear form on its dual which takes positive values on the
positive roots of $\fsl_3(\AA)$, and very large positive values on
those of $\fsl_3$. Then the positive roots of $\fg(\AA)$ are
those of $\fsl_3(\AA)$, those of $\fsl_3$, along with the weights
$\o_1+\m$, $\o_2-\m$ and $\o_1-\o_2-\m$, where $\m$ is a weight of 
$\fh_3(\AA)$. In particular, the highest root and the half-sum of the 
positive roots are 
$$\begin{array}{rcl}
\tilde{\a}_{\fg(\AA)} & = & \tilde{\a}_{\fsl_3}=\o_1+\o_2,\\
\rho_{\fg(\AA)}& = &
\rho_{\fsl_3(\AA)}+\rho_{\fsl_3}+\dim\fh_3(\AA)\o_1
  =  \rho_{\fsl_3(\AA)}+(3a+4)\o_1+\o_2.
\end{array}$$

We need to understand the distribution of the weights of the 
$\fsl_3(\AA)$-modules $\fh_3(\AA)$. They are as follows:


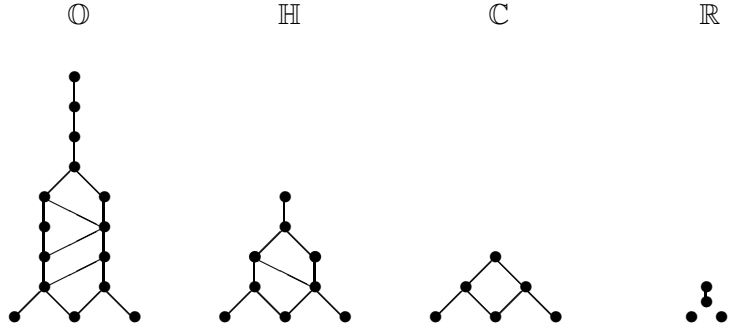
\setlength{\unitlength}{4mm}
\begin{picture}(30,11)(-8,0)
\put(0,10){${\mathbb O}$}
\put(0,5){$\bullet$}
\put(0,6){$\bullet$}
\put(0,7){$\bullet$}
\put(0,8){$\bullet$}
\put(0.2,5.2){\line(1,-1){1}}
\put(0.2,5.2){\line(-1,-1){1}}
\put(-1,4){$\bullet$}
\put(1,4){$\bullet$}
\put(0.2,5.2){\line(0,1){3}}
\put(-.8,1.2){\line(0,1){1}}
\put(-.8,2.2){\line(0,1){1}}
\put(-.8,3.2){\line(0,1){1}}
\put(1.2,1.2){\line(0,1){3}}
\put(-.8,4.2){\line(2,-1){2}}
\put(-.8,2.2){\line(2,1){2}}
\put(-.8,1.2){\line(2,1){2}}
\multiput(-1,1)(0,1){3}{$\bullet$}
\multiput(1,1)(0,1){3}{$\bullet$}
\multiput(-2,0)(2,0){3}{$\bullet$}
\put(-.8,1.2){\line(1,-1){1}}
\put(-.8,1.2){\line(-1,-1){1}}
\put(1.2,1.2){\line(1,-1){1}}
\put(1.2,1.2){\line(-1,-1){1}}

\put(7,10){${\mathbb H}$}
\put(7,4){$\bullet$}
\put(7,3){$\bullet$}
\put(7.2,3.2){\line(1,-1){1}}
\put(7.2,3.2){\line(-1,-1){1}}
\put(6,2){$\bullet$}
\put(8,2){$\bullet$}
\put(7.2,3.2){\line(0,1){1}}
\put(6.2,1.2){\line(0,1){1}}
\put(8.2,1.2){\line(0,1){1}}
\put(6.2,2.2){\line(2,-1){2}}
\multiput(6,1)(0,1){2}{$\bullet$}
\multiput(8,1)(0,1){2}{$\bullet$}
\multiput(5,0)(2,0){3}{$\bullet$}
\put(6.2,1.2){\line(1,-1){1}}
\put(6.2,1.2){\line(-1,-1){1}}
\put(8.2,1.2){\line(1,-1){1}}
\put(8.2,1.2){\line(-1,-1){1}}

\put(14,10){${\mathbb C}$}
\put(14,2){$\bullet$}
\put(14.2,2.2){\line(1,-1){1}}
\put(14.2,2.2){\line(-1,-1){1}}
\put(13,1){$\bullet$}
\put(15,1){$\bullet$}
\multiput(12,0)(2,0){3}{$\bullet$}
\put(13.2,1.2){\line(1,-1){1}}
\put(13.2,1.2){\line(-1,-1){1}}
\put(15.2,1.2){\line(1,-1){1}}
\put(15.2,1.2){\line(-1,-1){1}}

\put(21,10){${\mathbb R}$}
\put(21,1){$\bullet$}
\put(21,.5){$\bullet$}
\put(20.5,0){$\bullet$}
\put(21.5,0){$\bullet$}
\put(21.2,.6){\line(1,-1){.5}}
\put(21.2,.6){\line(-1,-1){.5}}
\put(21.2,.6){\line(0,1){.5}}
\end{picture}

\bigskip
The vertices of these diagrams indicate the weights with non-negative 
height (where the number $(\rho,\o)$ is the height of a weight
$\o$), while an edge indicates the action of a simple reflection
(the $\fh_3(\AA)$ are all minuscule modules, so that their sets of
weights are just the orbits of the highest ones). The complete diagram
is obtained by a symmetry along the line of height zero. 

The first three diagrams look very similar: there are three weights
of height zero, two weights on each height between $1$ and  
$\frac{a}{2}$, then one weight on each height up to $a$. 
This means that these three diagrams as being given by the superposition 
of intervals $[-a,a]$, $[-\frac{a}{2},\frac{a}{2}]$, plus a $0$. 
For $a=1$, this gives weights in height
$-1,-\frac{1}{2},0,\frac{1}{2},1$, with multiplicity two for the 
zero height: this is precisely our fourth diagram (where we have 
to use the normalization of \cite{bou} divided by two). 
 
\medskip We can analyse in a similar way our third series of models,
for which the weight distributions in the $\fsp_6(\AA)$-module
$\La^{\langle 3\rangle}\AA^6$ are again remarkably uniform. 
They are as follows: \bigskip


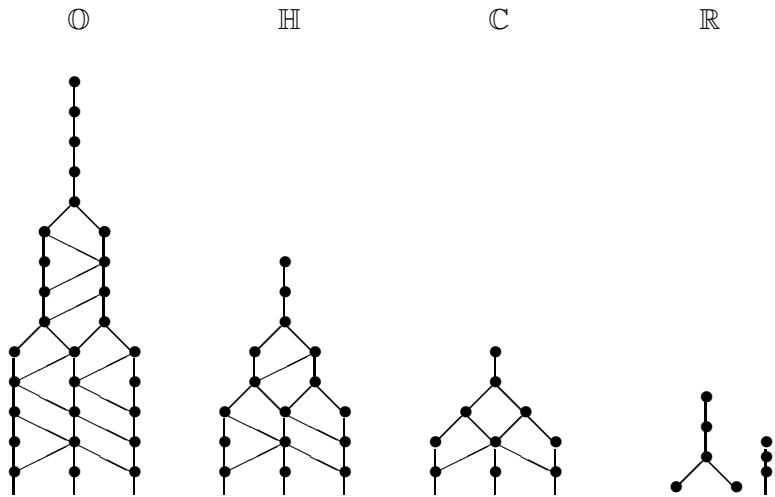
\setlength{\unitlength}{4mm}
\begin{picture}(30,16)(-8,0)
\put(0,15){${\mathbb O}$}
\put(0,10){$\bullet$}
\put(0,9){$\bullet$}
\put(0,11){$\bullet$}
\put(0,12){$\bullet$}
\put(0,13){$\bullet$}
\put(0.2,9.2){\line(0,1){4}}
\put(0.2,9.2){\line(1,-1){1}}
\put(0.2,9.2){\line(-1,-1){1}}
\put(-1,8){$\bullet$}
\put(1,8){$\bullet$}
\put(-.8,5.2){\line(0,1){3}}
\put(1.2,5.2){\line(0,1){3}}
\multiput(-1,5)(0,1){4}{$\bullet$}
\multiput(1,5)(0,1){4}{$\bullet$}
\put(-.8,8.2){\line(2,-1){2}}
\put(-.8,6.2){\line(2,1){2}}
\put(-.8,5.2){\line(2,1){2}}
\multiput(-2,0)(0,1){5}{$\bullet$}
\multiput(0,0)(0,1){5}{$\bullet$}
\multiput(2,0)(0,1){5}{$\bullet$}
\put(-.8,5.2){\line(1,-1){1}}
\put(-.8,5.2){\line(-1,-1){1}}
\put(1.2,5.2){\line(1,-1){1}}
\put(1.2,5.2){\line(-1,-1){1}}
\multiput(-1.8,-.5)(2,0){3}{\line(0,1){4.5}}
\multiput(-1.8,3.2)(2,0){2}{\line(2,1){2}}
\multiput(-1.8,3.2)(2,0){2}{\line(2,-1){2}}
\multiput(-1.8,2.2)(2,0){2}{\line(2,-1){2}}
\put(-1.8,.2){\line(2,1){2}}
\put(.2,1.2){\line(2,-1){2}}

\put(7,15){${\mathbb H}$}
\put(7,7){$\bullet$}
\put(7,6){$\bullet$}
\put(7,5){$\bullet$}
\put(7.2,5.2){\line(0,1){2}}
\put(7.2,5.2){\line(1,-1){1}}
\put(7.2,5.2){\line(-1,-1){1}}
\put(6.2,3.2){\line(0,1){1}}
\put(8.2,3.2){\line(0,1){1}}
\multiput(6,3)(0,1){2}{$\bullet$}
\multiput(8,3)(0,1){2}{$\bullet$}
\put(6.2,3.2){\line(2,1){2}}
\multiput(5,0)(0,1){3}{$\bullet$}
\multiput(7,0)(0,1){3}{$\bullet$}
\multiput(9,0)(0,1){3}{$\bullet$}
\put(6.2,3.2){\line(1,-1){1}}
\put(6.2,3.2){\line(-1,-1){1}}
\put(8.2,3.2){\line(1,-1){1}}
\put(8.2,3.2){\line(-1,-1){1}}
\multiput(5.2,-.5)(2,0){3}{\line(0,1){2.5}}
\multiput(5.2,2.2)(2,0){2}{\line(2,-1){2}}
\put(5.2,.2){\line(2,1){2}}
\put(7.2,1.2){\line(2,-1){2}}

\put(14,15){${\mathbb C}$}
\put(14,3){$\bullet$}
\put(14,4){$\bullet$}
\put(14.2,3.2){\line(0,1){1}}
\put(14.2,3.2){\line(1,-1){1}}
\put(14.2,3.2){\line(-1,-1){1}}
\multiput(13,2)(2,0){2}{$\bullet$}
\multiput(12,1)(2,0){3}{$\bullet$}
\multiput(12,0)(2,0){3}{$\bullet$}
\put(13.2,2.2){\line(1,-1){1}}
\put(13.2,2.2){\line(-1,-1){1}}
\put(15.2,2.2){\line(1,-1){1}}
\put(15.2,2.2){\line(-1,-1){1}}
\put(12.2,.2){\line(2,1){2}}
\put(14.2,1.2){\line(2,-1){2}}
\multiput(12.2,-.5)(2,0){3}{\line(0,1){1.5}}

\put(21,15){${\mathbb R}$}
\put(21,1.5){$\bullet$}
\put(21,.5){$\bullet$}
\put(21,2.5){$\bullet$}
\multiput(20,-.5)(2,0){2}{$\bullet$}
\put(21.2,.7){\line(1,-1){1}}
\put(21.2,.7){\line(-1,-1){1}}
\put(21.2,.7){\line(0,1){2}}
\multiput(23,0)(0,.5){3}{$\bullet$}
\put(23.2,-.5){\line(0,1){1.5}}
\end{picture}

\bigskip  
Again the fourth diagram is somewhat special: it splits into two orbits
of the Weyl group, the corresponding $\fsp_6$-module being non 
minuscule. Nevertheless, the first three diagrams are strikingly
similar: there are three strands of height from $\frac{1}{2}$ to 
$\frac{a+1}{2}$ , then two strands of length $\frac{a}{2}$, and 
a last strand of length $\frac{a}{2}+1$. Said otherwise, the
heights describe three intervals, namely
$[-\frac{3a+3}{2},\frac{3a+3}{2}]$, $[-\frac{2a+1}{2},\frac{2a+1}{2}]$
and $[-\frac{a+1}{2},\frac{a+1}{2}]$, and this, even for $a=1$. 
It is then very simple to apply the Weyl dimension formula to 
compute the dimension of $\fg(\AA)$ and its Cartan powers. 
A proof of proposition 1.1 stated in the introduction, 
valid for the entire exceptional series, follows.

\bigskip

\section{The general set up}

We say a collection of reductive Lie algebras $\fg(t)$  parametrized 
by $t$ and equipped with representations $(V_{\l_1}(t)\hd
V_{\l_p}(t))$ is a {\it series in strong the sense of Deligne} 
if there exists a  formula for $\tdim V_{m_1\l_1+\hdots+m_p\l_p}$
that is a rational function whose numerator and denominator are products
of linear functions of $t$. In this case, once one fixes $m_1\hd m_p$, 
the dimension formula  looks like  the Weyl dimension formula (see
below). We discuss other notions of series in \cite{lm4}. Note that
not all the dimension formulas of Vogel  are of this form
as algebraic 
extensions are required in his formulas,  see \cite{vog}.

\smallskip
How to construct such series?

\smallskip
One way would be to start with a fixed Lie algebra $\ff$, 
and consider $A$-graded Lie algebras $\fg$ (where $A$ is an
abelian group), containing $\ff$ as a component of $\fg_0$. 
If the grading comes from marking some nodes on the extended
Dynkin diagram of $\fg$, then $\ff$ will be given by a union 
of connected components  of the diagram obtained by removing  
the marked nodes. If one only marks one node, 
so one has a $\ZZ_2$-grading, then $\fg_0=\ff+\fh$ where $\fh$ 
is whatever else is left over after the nodes and components of 
$\ff$ are removed. In this   case,
$\fg_1=V\ot W(t)$ where $V$ (resp. $W$) is
the representation of $\ff$ (resp. $\fh(t)$) with highest weight the
sum of fundamental weights corresponding to nodes adjacent to the
marked node. For example, if one takes the node(s) next to the 
longest root, $\ff=\fsl_2$, and one can in particular recover the
last series of models of the exceptional Lie algebras in the preceding
section. If one takes the next node(s) over,
then $\ff=\fsl_3$ and one can recover the preceeding series. These two
gradings offer hope of universal formulas in the spirit of Vogel. 

\smallskip
In order to have a series in the strong sense
of Deligne, the Lie algebras   $\fh$ and representations
$U$  that remain must satisfy additional conditions explained below.

Write $\fg(t)= \ff + \fh(t) + W(t)$   
where  $\fg_0(t)=\ff+\fh(t)$ so $W(t)$ is a $\ff+\fh(t)$-module.
  We will   we need  that $W(t)=\Sigma_j V_j\ot U_j(t)$
where the $V_j$ are irreducible $\ff$-modules
all of the same dimension and  the $U_j(t)$ are irreducible
 $\fh(t)$-modules also all of the same dimension $u(t)$. We
will also need that $\trank \fg(t) = \trank \ff +\trank \fh(t)$
so we may chose Cartan subalgebras such that $\ft_{\fg}=\ft_f\op \ft_{\fh}$.
When there is no confusion, we supress the $t$.
The roots  of $\fg(t)$ are
\begin{itemize}
\item the roots of $\ff$,
\item the roots of $\fh$,
\item the weights $\mu+\nu$, with $\mu$ a weight of some  $V_j$ and $\nu$
a weight of $U_j$.
\end{itemize}
To get a set of positive roots we choose  linear forms $l$ and
$l_t$ on the root lattices, that are strictly positive on positive 
roots and heavily favor  the roots of $\ff$, so that the positive roots are:
\begin{itemize}
\item the positive roots of $\ff$,
\item the positive roots of $\fh$,
\item the weights $\mu+\nu$, with $\mu$ a weight of $V$ 
such that $l(\mu)>0$ and $\nu$
a weight of $U$. 
\end{itemize}
 
 We may write the half sum of the positive roots as
 $\rho_{\fg(t)}=\rho_f + \rho_{\fh} +u(t)\g$, where $\g$ is one half the
 sum of the positive weights of the $V_j$'s (positive in the sense that $l$ 
takes positive  values on them: we denote by $\Delta_+(V)$ the set of 
these weights).  
 We must reparametrize if necessary so that $u$ is a linear function of $t$.
 
  Now let $\o$ be a weight of $\fg$ suppported on 
$\g\in \ft_f$. This means that $\o$ is a weight of $\ff$ satisfying
 the integrality condition that $2(\o,\mu)/(\mu,\mu)\in\ZZ$
 for all $\mu\in\Delta_+(V)$.  (So in particular, $p$ above can at
 most be equal to the rank of $\ff$.) 
 
 We apply the Weyl dimension formula to $\o$. The contribution 
of the roots of $\fh$ to the prouct is trivial. 
The contribution of the roots of $\ff$ is

 $$\prod_{\a\in\Delta_+(\ff)}
\frac{(\rho_{\ff}+u(t)\g+\o,\a)}{(\rho_{\ff}+u(t)\g,\a)}
$$

The contribution of the other roots is more complicated, and  
to control this contribution   we  
add our most serious hypothesis: we require that when $t$ varies, 
the integers $(\rho_{\fg(t)},\m+\nu)$, 
for each set of values of $(\l_i,\m)$, not all zero, 
is the union  of intervals $[n_i(t)+1,m_i(t)]$,
where $n_i(t)$ and $m_i(t)$ are {\it linear} functions of $t$.
We allow  that for some values of $t$, $m_i(t)<n_i(t)$, 
which is to be interpreted as deleting the interval 
$[m_i(t)+1,n_i(t)]$. Then the contribution of such an interval 
to the Weyl dimension formula is:
 
$$
\frac{\binom{(\o,\m)+m_i}{(\o,\m)}}{\binom{(\o,\m)+n_i}{(\o,\m)}}.
$$

Putting these contributions together,
we see that we have a series in the strong sense of Deligne.

\bigskip

 \exam Here is a classical example. Let   
 $\fg(t)= \fso_{2t+4}$, $\ff=\fsl_2$, $\fh(t)=\fsl_2+\fso_{2t}$,
 $V=\CC^2$, $U=\CC^2\ot \CC^{2t}$.  
Let $\ff$ have root $\a$ and the $\fsl_2$ in $\fh(t)$ have root
$\b$. We use $\a_j$ to describe the roots of $\fso_{2l}$ and
sometimes the $\e_j$'s instead.
The positive roots  of $\fg(t)$ are
\begin{itemize}
\item $\a$,
\item $\b$, $\Delta_+(\fso_{2t})$
\item the weights $ \frac 12\a\pm\frac 12\b\pm\e_j$, $j=1\hd l$. 
\end{itemize}

We have $\rho_{\fg_f}=\frac12\a$, $\gamma= \frac 12 \a$ and
$\rho_{\fh}=\frac 12\b +(t-1)\e_1+(t-2)\e_2+\hdots +\e_{t-1}$. Thus,
taking inner products such that $(\e_i,\e_j)=\delta_{ij}$,
$(\a,\a)=(\b,\b)=2$, the pairings
$(\rho_{\fg(t)},\m+\nu)$ fill  the intervals $[1, 2t-1]$,
$[2, 2t]$, plus the isolated values $t$ and $t+1$. 
Thus applying our general formula we obtain
$$
\dim\fg\up k = 
\frac{(2k+2t+1)(k+t)(k+t+1)}{(2t+1)t(t+1)(k+1)}
\binom{k+2t-1}{k}\binom{k+2t}{k},
 $$
which is easy to obtain by directly applying the Weyl dimension 
formula. 
 \bigskip
 
 We conclude with an example of a two parameter series of Lie algebras
 in the strong sense of Deligne. 

\exam{The generalized third row of Freudenthal's magic chart.} 
With the notations of section 4, we have
$$
\fg(r,a)=\fg_r (\AA,\HH) = \ft_r(\AA)\op \fsl_2^{\times r}\op 
\Sigma_{i<j}U_i\ot U_j\ot\AA_{ij}.
$$
We have $\fg(r,a)= \fsl_2^{xr}$ for $a=0$, $fc_r$ when $a=1$,
$\fa_{2r-1}$ when $a=2$, $\fd_{2r}$ when $a=4$ and $\fe_7$ when
$a=8$ and $r=3$.

With our conventions, the positive roots of $\fg(t)$ are:
\begin{itemize}
\item the positive roots $\a_i$, $1\le i\le r$,  of $\fsl_2^{\times r}$,
\item   the positive roots of $\ft_r(\AA)$,
\item the weights $\o_i-\o_j+\mu_{ij}$, $i<j$.
\end{itemize}
Write the half sum of positive roots as $\rho =\rho_{\ft_r(\HH)}+
\rho_{\ft_r(\AA)} + a\g$
with $2\g=(r-1)\a_1+(r-2)\a_2+\cdots +\a_{r-1}$. 
Applying our method once again, we obtain the three parameters
formula
$$\dim\fg_r(\HH,\AA)^{(k)}=
\frac{2k+a(r-1)+1}{a(r-1)+1}
\frac{\binom{k+\frac{ar}{2}-1}{k}\binom{k+ar-a}{k}
\binom{k+\frac{ar-a}{2}}{k}\binom{k+ar+1-\frac{3a}{2}}{k}}
{\binom{k+\frac{a}{2}-1}{k}
\binom{k+\frac{ar}{2}+1-a}{k}\binom{k+\frac{ar-a}{2}}{k}}.
$$

One can derive similar formulas for all the preferred representations
in the generalized second and third rows.

\end{document}